\documentclass[11pt]{article}
\usepackage[cm]{fullpage}\usepackage{latexsym}\usepackage{amsfonts}\usepackage{amsmath}\usepackage{amssymb}\usepackage{amscd}
 \newcommand{\be}{\begin{equation}}
       \newcommand{\ee}{\end{equation}}
       \newcommand{\ba}{\begin{eqnarray}}
        \newcommand{\ea}{\end{eqnarray}}
 \newcommand{\ban}{\begin{eqnarray*}}
 \newcommand{\ean}{\end{eqnarray*}}

 \newcommand{\ra}{\rightarrow}

 \newcommand{\vol}{\mathrm{Vol}}

    \newcommand{\qed}{\hspace*{\fill}\rule{3mm}{3mm}\quad \vspace{.2cm}}
  \newcommand{\Pf}{\noindent {\bf Proof:} }

 \newcommand{\sect}[1]{\section{#1} \setcounter{equation}{0}}

\newtheorem{theorem}{Theorem}[section]

\newtheorem{remark}[theorem]{Remark}

\newtheorem{proposition}[theorem]{Proposition}
\newtheorem{question}[theorem]{Question}

\begin{document}

\title{\normalsize{ON VOLUME GROWTH OF GRADIENT STEADY RICCI SOLITONS}}
\author{\scriptsize{GUOFANG WEI, PENG WU}}
\date{}
\maketitle
\begin{abstract}
In this paper we study volume growth of gradient steady Ricci
solitons. We show that if the potential function satisfies a uniform
condition, then the soliton has at most Euclidean volume growth.
\end{abstract}

\sect{Introduction}

$(M^n, g)$ is a gradient  Ricci soliton if there is smooth function
$f: M \ra \mathbb R$ and constant $\lambda \in \mathbb R$ such that
\begin{equation}
\text{Ric}+\text{Hess} f = \lambda g.
\end{equation}
$f$ will be referred as the potential function. The soliton is
called shrinking, steady, expanding when $\lambda >0, \lambda = 0,
\lambda <0$ respectively.

Ricci solitons are self-similar solutions of the Ricci flow, and
play an important role in the study of singularity formation. They
are also natural extensions of Einstein manifolds, and special cases
of smooth metric measure spaces.

Volume growth of gradient Ricci solitons is of particular interest
to mathematicians. Estimate of the potential functions plays an
important role in the study of volume growth. In
\cite{Hamilton1995}, Hamilton proved the following identity for
gradient Ricci solitons,
\[
R+|\nabla f|^2 -2\lambda f = \Lambda,
\]
where $\Lambda$ is a constant, and $R$ is the scalar curvature.

For gradient shrinking Ricci solitons, the answer is complete.
Perelman \cite{Perelman} and Cao-Zhou \cite{Cao-Zhou2010} proved
that $f$ always grows quadratically. Cao-Zhou \cite{Cao-Zhou2010}
further proved that any gradient Ricci shrinking soliton has at most
Euclidean volume growth. Recently, Munteanu-Wang \cite{MW} proved
that  any gradient Ricci shrinking soliton has at least linear
volume growth.

For gradient steady Ricci solitons, B.-L. Chen \cite{Chen2009}
proved that $R\ge 0$. Hence $\Lambda \ge 0$ and equal to zero if and
only if $f$ is constant and $(M, g)$ is Ricci flat. When $\Lambda
>0$ we can assume $\Lambda = 1$ after scaling, i.e. \be R+|\nabla f|^2 =1.
\label{steady-s} \ee Combine with the trace of steady Ricci soliton
equation $R+\Delta  f=0$, we have \be \Delta f - |\nabla f|^2 = -1.
\label{Delta-f} \ee Therefore $f$ has no local minimum.
(\ref{steady-s}) and that $R \ge 0$ also gives $|\nabla f| \le 1$.
Namely $f$ decays at most linearly.

Cao-Chen \cite{Cao-Chen} proved that $f$ decays linearly when Ricci
curvature is positive and $R$ attains its maximal at some point.
However the simple example of $\mathbb R^2$ with the canonical
metric $g_0$ and $f(x) = x_1$ shows that this is not the case: $f$
is constant along $x_2$ direction. Note that Riemannian product of
any two steady gradient Ricci solitons is still a steady gradient
Ricci soliton. Hence a steady gradient Ricci soliton multiply with a
trivial one ($f$ is a constant) will have constant direction. Though
one can make product of two shrinking ones, but all trivial
shrinking ones are compact. So it will not give a constant direction
by taking a product. Munteanu-Sesum \cite{MS} and the second author
\cite{Wu} independently showed that the infimum of $f$ does decay
linearly. In fact
\begin{equation}
\begin{split}
-r &\leq \inf_{y\in\partial B(x,r)}f(y)-f(x) \leq -r
+\sqrt{2n}(\sqrt{r}+1),\hspace{1cm} r\gg 1. \label{sqrt-r}
\end{split}
\end{equation}
In particular, $\liminf_{y\ra \infty} R(y) = 0$, see also
\cite{Fernadez-Garcia2011, Chow-Lu-Yang2011}.

We note that among all known examples of steady gradient Ricci
solitons, the infimum of $f$ is like $-r + O(\ln r)$. See the survey
article \cite{Cao2010} for a list of examples.  One naturally asks
if one can improve the second order term in (\ref{sqrt-r}) from
$\sqrt{r}$ to $\ln r$. We shows this is indeed the case for a large
class of steady gradient Ricci solitons.  To study the second order
term, write the potential function in polar coordinate:
\begin{equation*}
\begin{split}
f(r,\theta)=-r+\phi(r,\theta)
\end{split}
\end{equation*}
where $r(\cdot)=d(x,\cdot)$ for some $x\in M^n$, $\theta\in
S^{n-1}$. Without loss of generality we assume $\phi(0,\theta)=0$ by
adding a constant to $f$. Since $f(r) \ge -r$ and $|\nabla f|\leq
1$, $\phi(r) \ge 0$ and $\phi(r,\theta)$ is nondecreasing in $r$ for
any fixed $\theta$. We show that the estimate (\ref{sqrt-r}) can be
improved to $\ln r$, if $\phi$ at one direction is comparable to
minimum of $\phi$ among all spherical directions for all $r$ large.
Namely,

\begin{theorem} Let $(M^n,\, g,\, f)$ be a complete
gradient steady Ricci soliton satisfying (\ref{steady-s}). Assume
that there exist $\theta_0\in S^{n-1}$, and constants $C_1 \ge 0,
C_2 \ge 0$ such that
\begin{equation}
\begin{split}
\int_0^r  \left( \phi(r,\theta_0)-\phi(t,\theta_0) \right) dt &\leq
C_1\min_{\theta\in S^{n-1}}\int_0^r \left(
\phi(r,\theta)-\phi(t,\theta) \right) dt +C_2 r  \label{min-dir}
\end{split}
\end{equation}
for sufficiently large $r$. Then for any $x\in M^n$, there exist
constants $C\ge 0$, $r_0>0$ such that for  $r\geq r_0$,
\begin{equation}
\begin{split}
-r \leq  \inf_{y\in\partial B(x,r)}f(y)-f(x) &\leq -r + \left(\frac
n2 C_1 + C_2\right) \ln r + C. \label{inf-f}
\end{split}
\end{equation}
\end{theorem}

All known examples of gradient steady Ricci solitons satisfy the
condition (\ref{min-dir}). We suspect that the estimate
(\ref{inf-f}) holds for all gradient steady Ricci solitons.

In \cite{MS}, Munteanu and Sesum proved that any gradient steady
Ricci soliton has at least linear volume growth, and at most growth
rate of $e^{\sqrt{r}}$. We show that if the potential function
satisfies a uniform condition in the spherical directions, then the
gradient steady Ricci soliton has at most Euclidean volume growth.

\begin{theorem}
Let $(M^n,\, g,\, f)$ be a complete gradient steady Ricci soliton
satisfying (\ref{steady-s}). Assume that there exist constants
$C_1,C_2 \ge 0$ such that
\begin{equation}
\begin{split}
\max_{\theta\in S^{n-1}}\int_0^r \phi(r,\theta)-\phi(t,\theta) dt
&\leq C_1\min_{\theta\in S^{n-1}}\int_0^r
\phi(r,\theta)-\phi(t,\theta) dt +C_2 r  \label{f-uniform}
\end{split}
\end{equation}
for sufficiently large $r$. Then for any $x\in M^n$, there exists
$r_0>0$, for any $r\geq r_0$,
\begin{equation}
\begin{split}
-r \leq  f(y)-f(x) &\leq -r + C\ln r.
\end{split}
\end{equation}
for any $y\in \partial B(x,r)$. Moreover, the soliton has at most
Euclidean volume growth, i.e. for any $x\in M^n$, there exists
$r_0>0$, for any $r\geq r_0$,
\begin{equation*}
\begin{split}
\text{Vol}(B(x,r)) &\leq Cr^n.
\end{split}
\end{equation*}
If in addition $\phi(r)\ge \delta\ln r$ for large $r$, then
\begin{equation*}
\begin{split}
\text{Vol}(B(x,r)) &\leq Cr^{n-\delta}.
\end{split}
\end{equation*}
\end{theorem}

\begin{remark}
(1). If $\phi$ increases uniformly along all spherical directions,
i.e. \\
$\max_{\theta} \frac{\partial \phi}{\partial r}\leq C\min_{\theta}
\frac{\partial \phi}{\partial r}$, where $\theta\in S^{n-1}$, then
$\phi$ satisfies (1.7) with $C_1=C,\ C_2=0$.

(2). Theorem 1.2 can be considered as an analogue of volume growth
theorem for gradient shrinking Ricci solitons of
\cite{Cao-Zhou2010}. For gradient shrinking Ricci solitons, the
potential function automatically satisfies a uniform condition
\cite{Cao-Zhou2010}; for gradient steady Ricci solitons, we need to
impose a uniform condition.

(3). If the soliton is rectifiable (see \cite{PW}), i.e. $f$ is the
distance function from a set, then $\phi$ satisfies
(\ref{f-uniform}) with $C_1 = 1$ if the set is bounded (this is the
case with all nonproduct examples).
\end{remark}

To prove the results, the following estimate for $\phi$ which holds
for all gradient steady Ricci solitons is the key.
\begin{proposition} \label{key-est}
Let $(M^n,\, g,\, f)$ be a complete gradient steady Ricci soliton
satisfying (\ref{steady-s}). Then \be \min_{y\in \partial B(x,r)}
 \int_0^{r} (\phi(y) - \phi(t)) dt \le \frac{n}{2}(r+
\sqrt{r}) + o(\frac 1r). \label{pre-min} \ee
\end{proposition}

This estimate improves the estimate in \cite{Wu}. In the next
section we derive a volume comparison for the solitons by adapting
the volume comparison for smooth metric measure in \cite{WW}. Then
we prove Proposition~\ref{key-est} by combining with equation
(\ref{Delta-f}). In Section 3 we prove the main theorems using this
estimate and an ODE.

\textbf{Acknowledgments.} The second author would like to thank
Professor Thomas Sideris for helpful discussions.

\section{The Preliminary Estimate}

In this section we prove Proposition~\ref{key-est} by applying an
weighted volume comparison argument for smooth metric measure spaces
as in \cite{WW, Wu}.

Recall a smooth metric measure space is a triple $(M^n,\ g,\
e^{-f}\text{dvol}_g)$, where $(M^n,\ g)$ is a smooth Riemannian
manifold, and $f:\ M^n\rightarrow \mathbb{R}$ is a smooth function.
Write the volume element in polar coordinate $dvol = J(r,\theta) dr
d\theta$. Define the weighted volume element $J_f(r, \theta) =
e^{-f} J(r, \theta)$, and weighted volume $vol_f B(x,r) =
\int_{B(x,r)} e^{-f} dvol$.

Wei and Wylie \cite{WW} obtained the following $f$-volume comparison
theorem for smooth metric measure spaces,

\begin{theorem} ($f$-volume comparison) \\
Suppose $(M^n,g,e^{-f}dvol)$ is a smooth metric measure space with
$Ric_f\geq (n-1)H$. Fix $x\in M$. If $|f|\leq \Lambda$, then for
$R\geq r>0$ ($R\leq \pi\slash 4\sqrt{H}$ if $H>0$),
\begin{equation*}
\begin{split}
\frac{V_f(B_R(x))}{V_f(B_r(x))} &\leq
\frac{V_H^{n+4\Lambda}(B_R)}{V_H^{n+4\Lambda}(B_r)}.
\end{split}
\end{equation*}
Where $V^n_H(B_r)$ is the volume of the ball of radius $r$ in
$M^n_H$, the simply connected model space of dimension $n$ with
constant sectional curvature $H$.
\end{theorem}

One observes that the dimension of the model space in the volume
comparison depends on the potential function $f$. A further
investigation of the dimension will lead to Proposition
\ref{key-est}.

Denote $m_f = (\ln J_f)'$, the $f$-mean curvature. For $0 <
r_1 \le r_2$, let $A(x,r_1, r_2) = \{y| r_1 \le d(x,y) \le r_2\}$ be
the annulus,  and $$a = \min_{y\in A(x,r_1, r_2)} \frac{2}{r(y)^2}
\int_0^{r(y)} (\phi(y) - \phi(t)) dt.$$ Clearly $a \ge 0$. By
(\ref{sqrt-r}), we have $a\leq \frac{C}{\sqrt{r_1}}$ for $r_1 \gg
1$. For the rest we assume $r_1 \gg 1$ and therefore we can assume
$a<1$.
\begin{proposition}
For gradient steady Ricci soliton, we have \be  m_f(r,\theta) \leq
\frac{n-1}{r}+ 1-\frac{2}{r^2}\int_0^r
[\phi(r,\theta)-\phi(t,\theta)]dt \le \frac{n-1}{r}+ 1.
\label{m-f-est}
\end{equation}
and \begin{equation} \frac{\vol_f (\partial B(x, r_2))}{\vol_f
(A(x,r_1,r_2))} \le \frac{\frac{n}{r_2} + 1-a}{1-
(\frac{r_1}{r_2})^{n+(1-a)r_2}}.
 \label{vol-le}
\end{equation}
\end{proposition}

\Pf  For smooth metric space  $(M^n,g,f)$ with $\mbox{Ric}_f \ge 0$,
recall the following estimate for $m_f$  from
 \cite[(3.19)]{WW},
\begin{equation*}
\begin{split}
m_f(r,\theta) &\leq \frac{n-1}{r}+ \frac{2}{r^2} \int_0^r (f(t)-
f(r)) dt.
\end{split}
\end{equation*}
Plug in $f = -r + \phi$ gives (\ref{m-f-est}).

Now let
\[ \overline{m} (r) = \left\{
\begin{array}{ll}
\frac{n-1}{r} + 1 & r \le r_1 \\
\frac{n-1 + (1-a)r_2}{r} & r_1 < r \le r_2
\end{array} \right.,
\]
then \begin{equation} m_f(r) \le \overline{m} (r)  \ \ \mbox{for} \
0<r \le r_2. \label{mean-c}  \end{equation}

Let $\overline{A}(r)= e^{\int_0^r \overline{m}(t)dt}$ and
$\overline{V}(r_0, r) = \int_{r_0}^r \overline{A}(t)dt$. From the
mean curvature relation (\ref{mean-c}), we have $\left(
\frac{A_f}{\overline{A}}\right)' \le 0$, therefore
\[ \frac{\vol_f (\partial B(x, r_2))}{\vol_f
(A(x,r_1,r_2))} \le \frac{\overline{A}(r_2)}{\overline{V}(r_1,r_2)},
\]

We compute \ban \frac{\overline{A}(r_2)}{\overline{V}(r_1,r_2)} & =
& \frac{e^{\int_{0}^{r_2} \overline{m}(t) dt}}{\int_{r_1}^{r_2}
e^{\int_{0}^{s} \overline{m}(t) dt} ds} = \frac{e^{\int_{r_1}^{r_2}
\overline{m}(t) dt}}{\int_{r_1}^{r_2}
e^{\int_{r_1}^{s} \overline{m}(t) dt} ds} \\
& = & \frac{(r_2/r_1)^{n-1+(1-a)r_2}}{\int_{r_1}^{r_2}(s/r_1)^{n-1+(1-a)r_2} \, ds} \\
& = & \frac{\frac{n}{r_2} + 1-a}{1- (\frac{r_1}{r_2})^{n+(1-a)r_2}}.
\ean
This gives (\ref{vol-le}). \qed

\noindent{\bf Proof of Proposition~\ref{key-est}}.

Integrating (\ref{Delta-f}) and using
 $|\nabla f| \le 1$, we have, for any $x \in M$, \ban
\int_{B(x, r)} 1 \cdot e^{-f} dvol & = & -  \int_{B(x, r)} (\Delta f - |\nabla f|^2)  \cdot e^{-f} dvol \\
&= & - \int_{\partial B(x, r)} \frac{\partial f}{\partial n} e^{-f}
dvol \le   \int_{\partial B(x, r)}  e^{-f} dvol. \ean Therefore \be
\frac{\vol_f (\partial B(x,r))}{\vol_f (B(x,r))} \ge 1.
\label{vol-ge}\ee

Combining (\ref{vol-le}) and
 (\ref{vol-ge}) we have
 \[a \le \frac{n}{r_2} + (\frac{r_1}{r_2})^{n+(1-a)r_2}.\]
 Let $r_1 = r, r_2 = r+\sqrt{r}$, then $r_1/r_2 =
 (1+\frac{1}{\sqrt{r}})^{-1}$. When $r$ is large,
 \[ (1+\frac{1}{\sqrt{r}})^{-(n+(1-a)(r+\sqrt{r}))} = O(e^{-(1-a)\sqrt{r}}).\]
 Therefore, for all $r$ large enough,
\[
a = \min_{y\in A(x,r, r+\sqrt{r})} \frac{2}{r(y)^2} \int_0^{r(y)}
(\phi(y) - \phi(t)) dt  \le \frac{n}{r+ \sqrt{r}} +
O(e^{-(1-a)\sqrt{r}}).
\]
Suppose the minimum in above is attained at $y_0 = (r_0,\theta_1)$
with $r \le r_0 \le r+\sqrt{r}$. Then \ban \min_{y\in
\partial B(x,r)} \int_0^{r} (\phi(y) - \phi(t)) dt  & \le & \int_0^{r_0}
(\phi(y_0) - \phi(t)) dt  \\
& \le &  \frac{r_0^2}{2} \left( \frac{n}{r+ \sqrt{r}} +
O(e^{-(1-a)\sqrt{r}}) \right) \\
& \le & \frac{n}{2}(r+ \sqrt{r}) + o(\frac 1r). \ean
 \qed

\section{Proof of Main Results}

\noindent {\bf Proof of Theorem 1.1}. From (\ref{pre-min}) and
(\ref{min-dir}),  we have
\begin{equation}
\int_0^r [\phi(r,\theta_0)-\phi(t,\theta_0)]dt \leq
\frac{nC_1}{2}(r+ \sqrt{r})+C_2r + o(\frac 1r). \label{phi-est}
\end{equation}
For simplicity when there is no confusion we omit $\theta_0$ in the
function. Let $\Phi(r)=\int_0^r \phi(t)dt$, then the above
inequality can be written as,
\begin{equation}
\begin{split}
\Phi'(r)-\frac{1}{r}\Phi(r) &\leq \frac{nC_1}{2}+C_2 +
O(\frac{1}{\sqrt{r}}).  \label{Phi}
\end{split}
\end{equation}
Multiply by the integrating factor $\frac{1}{r}$ and integrate from
some fixed $t_0 \gg 1$ to $r$, we get
\begin{equation*}
\begin{split}
\frac{\Phi(r)}{r} &\leq \left(\frac{nC_1}{2}+C_2\right)\ln r +C_3.
\end{split}
\end{equation*}
So we have \ban \phi(r,\theta_0) &= & \Phi'(r,\theta_0)\leq
\frac{\Phi(r,\theta_0)}{r}+\frac{nC_1}{2}+C_2 +
O(\frac{1}{\sqrt{r}}) \\
&  \leq & \left(\frac{nC_1}{2}+C_2\right)\ln r  + C_4\\
f(r,\theta_0) &= & -r+\phi(r,\theta_0) \leq -r +
\left(\frac{nC_1}{2}+C_2\right)\ln r  + C_4. \ean This gives
(\ref{inf-f}). \qed

\

\noindent {\bf Proof of Theorem 1.2} From (\ref{pre-min}) and
(\ref{f-uniform}), we have
\begin{equation*}
\begin{split}
\int_0^r [\phi(r,\theta)-\phi(t,\theta)]dt &\leq  \frac{nC_1}{2}(r+
\sqrt{r})+C_2r + o(\frac 1r),
\end{split}
\end{equation*}
for all $\theta\in S^{n-1}$. Therefore (\ref{inf-f}) holds for all
$y$. Namely, for all $y\in\partial B_r(x)$,
\begin{equation*}
\begin{split}
-r\leq f(y)-f(x) &\leq -r + \left(\frac{nC_1}{2}+C_2\right)\ln r +
C_4.
\end{split}
\end{equation*}

By (\ref{m-f-est}), for all $r>0$,
\begin{equation*}
\begin{split}
\frac{\partial}{\partial r}\ln J &=m_f(r)+\langle\nabla f,\nabla r\rangle\\
&\leq \frac{n-1}{r}+1-\frac{2}{r}\phi (r) +  \frac{2}{r^2}\int_0^r
\phi(t) dt +\langle\nabla f,\nabla r\rangle.
\end{split}
\end{equation*}
Integrate from 1 to r and do integration by part for the double
integral, we get \ba \ln J(r)-\ln J(1) & \leq& (n-1)\ln
r+(r-1)-\int_1^r \frac 2s \phi(s) ds \nonumber \\
& & +\left( - \frac{2}{s} \int_0^s \phi(t) dt \right) \left|_1^r \right. + \int_1^r \frac 2s \phi(s) ds  +f(r)-f(1) \nonumber \\
& =& (n-1)\ln r  +\phi(r) - \frac{2}{r}\int_0^r \phi(t)dt +
2\int_0^1 \phi(t) dt-f(1) \label{vol}
 \\
 & = & (n-1)\ln r  - \phi(r) +2\left(\phi(r) - \frac{1}{r}\int_0^r \phi(t)dt \right) + 2\int_0^1 \phi(t) dt-f(1). \nonumber
\ea
 Using (\ref{Phi}) we have, for $r$ large,
\[
\ln J(r) \le   (n-1)\ln r  - \phi(r) + C \le (n-1)\ln r  + C.
\]
Hence
\begin{equation*}
\begin{split}
J(r) &\leq e^Ce^{(n-1)\ln r}=e^Cr^{n-1},
\end{split}
\end{equation*}
and  the volume of geodesic ball centering at $x$ satisfies
\begin{equation*}
\begin{split}
\text{Vol}(B(x,r))&\leq C' r^{n}.
\end{split}
\end{equation*}

If further $\phi(s) \ge \delta \ln s$, then we have
\begin{equation*}
\begin{split}
J(r) &\leq  Cr^{n-1} \exp(-\phi(r))\leq Cr^{n-\delta-1},
\end{split}
\end{equation*}
therefore the volume growth is strictly less than Euclidean volume
growth,
\begin{equation}
\begin{split}
\text{Vol}(B(x,r)) &\leq Cr^{n-\delta}.
\end{split}
\end{equation}
\qed

For general gradient steady Ricci solitons, the estimate of
potential function can be reduced to the following
\begin{question}
Suppose $\phi$ is nondecreasing along any minimal geodesic starting
from $x$. Assume that for $r$ sufficiently large,
$\inf_{y\in\partial B(x,r)}\phi(y)\leq C\sqrt{r}$ and
\begin{equation*}
\begin{split}
\inf_{y\in\partial B(x,r)}\int_1^r[\phi(y)-\phi(\gamma_y(t)]dt &\leq
\frac{nr}{2}.
\end{split}
\end{equation*}
Does the following hold?
\begin{equation*}
\begin{split}
\inf_{y\in\partial B(x,r)}\phi(y)\leq C\ln r ?
\end{split}
\end{equation*}
\end{question}

\

\begin{remark}
From (\ref{vol}), we see that if
\begin{equation*}
\begin{split}
-r \leq  f(y)-f(x) &\leq -r + C\ln r.
\end{split}
\end{equation*}
for $y\in \partial B(x,r)$, then for any $x\in M^n$, there exists
$r_0>0$ such that for any $r\geq r_0$,
\begin{equation*}
\begin{split}
\text{Vol}(B(x,r)) &\leq C'r^{n+C}.
\end{split}
\end{equation*}
\end{remark}


\

Guofang Wei, Math Department, UCSB, Santa Barbara, CA 93106

wei@math.ucsb.edu

\

Peng Wu, Math Department, UCSB, Santa Barbara, CA 93106

Current address: Math Department, Cornell University, Ithaca, NY
14853.

wupenguin@math.ucsb.edu


\begin{thebibliography}{}
\renewcommand{\baselinestretch}{0}
\footnotesize

\bibitem{Chow-Lu-Yang2011} Bennett Chow, Peng Lu, and Bo Yang, \emph{ Lower bounds for the scalar
curvatures of noncompact gradient Ricci solitons.} C. R. Math. Acad.
Sci. Paris 349 (2011), no. 23-24, 1265--1267,

\bibitem{Cao2010}Huai-Dong Cao, \emph{Recent progress on Ricci solitons.} Recent advances in geometric analysis,
1--38, Adv. Lect. Math. (ALM), 11, Int. Press, Somerville, MA, 2010.

\bibitem{Cao-Chen} Huai-Dong Cao and Qiang Chen, \emph{On Locally Conformally Flat Gradient Steady Ricci Solitons},
arXiv:math.DG/0909.2833, to appear in Trans. Amer. Math. Soci.

\bibitem{Cao-Zhou2010} Huai-Dong Cao and Detang Zhou, \emph{On Complete Gradient Shrinking Ricci Solitons},
J. Differential Geom. \textbf{85} (2010), 175--186.

\bibitem{Chen2009} Bing-Long Chen, \emph{Strong Uniqueness of the Ricci Flow}, J. Differential
Geom. \textbf{82} (2009), 363--382.

\bibitem{Fernadez-Garcia2011} Manuel Fern\'andez-L\'opez and Eduardo Garc\'ia-R\'io, \emph{Maximum
principles and gradient Ricci solitons}, J. Differential Equations
251 (2011), no. 1, 73--81.

\bibitem{Hamilton1995} Richard Hamilton, \emph{The formation of singularities in the Ricci flow}, Surveys in Differential
Geometry (Cambridge, MA, 1993), 2, 7-136, International Press,
Combridge, MA,1995.

\bibitem{MS} Ovidiu Munteanu and Natasa Sesum, \emph{On Gradient Ricci Solitons},
arXiv:math.DG/0910.1105, to appear in J. Geom. Anal.

\bibitem{MW} Ovidiu Munteanu and Jiaping Wang, \emph{Analysis of the weighted Laplacian and applications to Ricci solitons},
arXiv:math.DG/1112.3027, to appear in Comm. Anal. Geom.

\bibitem{Perelman} Grisha Perelman, \emph{Ricci flow with surgery on three manifolds}, arXiv:math.DG/0303109.

\bibitem{PW} Peter Petersen and William Wylie, \emph{On gradient Ricci solitons with symmetry}, Proc. Amer. Math.
Soc. \textbf{137} (2009), 2085--2092.

\bibitem{WW} Guofang Wei and William Wylie, \emph{Comparison Geometry for the Bakry-Emery Ricci Tensor}, J.
Differential Geometry \textbf{83} (2009), 377--405.

\bibitem{Wu} Peng Wu, \emph{On the Potential Function of Gradient Steady Ricci Solitons}, J. Geom. Anal. 2011, DOI 10.1007/s12220-011-9243-7.


\end{thebibliography}
\end{document}